\documentclass{article}
\usepackage[utf8]{inputenc}
\usepackage[english]{babel}
\usepackage{amsmath,amsthm,amsfonts}
\usepackage{amssymb}
\usepackage{graphicx} 
\usepackage{bm}
\usepackage{tikz}
\usepackage{multirow}
\usepackage{color}
\usepackage{placeins}
\usepackage{graphicx}
\usepackage{cite}

\usepackage{multicol}
\usepackage{wrapfig,lipsum,booktabs}

\graphicspath{ {./figs/} }

\begin{document}

\title{Non-Local Multi-Continuum method (NLMC) for Darcy-Forchheimer flow  in fractured media}

%

\author{
Denis Spiridonov
\thanks{Laboratory of Computational Technologies for Modeling Multiphysical and Multiscale Permafrost Processes, North-Eastern Federal University, 677000 Yakutsk,Republic of Sakha (Yakutia), Russia. 
Email: {\tt d.stalnov@mail.ru}.}
\and
Maria Vasilyeva \thanks{Department of Mathematics and Statistics, Texas A\&M University, Corpus Christi, Texas, USA. 
Email: {\tt maria.vasilyeva@tamucc.edu}.}
}




\maketitle

\begin{abstract}

This work presents the application of the non-local multicontinuum method (NLMC)  for the Darcy-Forchheimer model in fractured media. The mathematical model describes a nonlinear flow in fractured porous media with a high inertial effect and flow speed. The space approximation is constructed on the sufficiently fine grid using a finite volume method (FVM) with an embedded fracture model (EFM) to approximate lower dimensional fractures. A non-local model reduction approach is presented based on localization and constraint energy minimization. The multiscale basis functions are constructed in oversampled local domains to consider the flow effects from neighboring local domains. Numerical results are presented for a two-dimensional formulation with two test cases of heterogeneity. The influence of model nonlinearity on the multiscale method accuracy is investigated. The numerical results show that the non-local multicontinuum method provides highly accurate results for  Darcy-Forchheimer flow in fractured media.

\end{abstract}

\section{Introduction}

Fluid flow in porous media is an essential element in understanding oil and gas production processes, as well as in reservoir hydrology, environmental protection, and many other \cite{cescon2020filtration, maitland2000oil, schofield2000gas}.
Properly describing fluid properties and their movement in the reservoir is necessary for the model construction process to solve the problems of designing field development systems. This task requires accurately describing fluid and gas flow under actual reservoir conditions. Mathematical modeling is often used to select the best option for oil and gas field development. Modeling helps test various hydrocarbon production technologies, find the best options for production well placement schemes, and determine changes in physical oil and gas parameters \cite{d2007mathematical, he2009progress, reite2006mathematical}. 

The fundamental law of fluid flow in porous media is Darcy's law; it expresses the dependence of fluid filtration rate on the pressure gradient. Many scientific works are devoted to checking and investigating the limits of applicability of Darcy's law \cite{gavrilieva2017flow, vasil2018numerical, nikiforov2019numerical}. In the case when the filtration rate is relatively high, inertial effects cannot be ignored, and the Darcy-Forchheimer model should be used \cite{lee1997modeling, khan2020simulation, chu2021stability}. Moreover, a high filtration rate can occur in highly heterogeneous and fractured media. In the fractures continuum, the flow velocity is much higher than in porous media \cite{vasilyeva2022efficient}.
In this paper, we consider the multiscale method for the nonlinear flow model in a fractured porous medium based on the Darcy-Forchheimer model \cite{girault2008numerical}.

Standard modeling approaches imply well-known approximation methods on a fine grid. The finite element or finite volume method is well-proven in mathematical modeling \cite{bathe2007finite, dhatt2012finite, jagota2013finite}. For modeling the Darcy-Forchheimer model, the mixed finite element method is commonly used to preserve the mass conservation on the discrete level \cite{salas2013analysis, huang2018multigrid, xu2017multipoint}. This paper uses the finite volume method to solve a problem on the fine grid \cite{eymard2000finite, moukalled2016finite, barth2003finite, karimi2009detailed}. For the problem in fractured media, we should carefully choose the method for the approximation. The most straightforward way to approximate fractures based on the explicit fracture representation on the grid with the application of the discrete fracture model (DFM) for their approximation  \cite{jaffre2011discrete, vasil2017numerical, garipov2016discrete}. However, the DFM approach requires the construction of a very detailed grid with a vast number of cells for the case of extensive fracture distribution, which leads to substantial computational costs. In this work, we use an embedded fractured model (EFM) that allows constructing mesh for fracture networks independently of porous media mesh \cite{spiridonov2017simulation, shakiba2015using, tyrylgin2019embedded}. However, a detailed fine grid is still required for a heterogeneous porous media with high contrast coefficients \cite{ vasilyeva2020learning, vasilyeva2021machine, vasilyeva2021preconditioning}.

One way to solve this problem is to solve the problem on a coarse grid \cite{grigoriev2019effective, stepanov2023prediction, tyrylgin2019numerical}. In this study, we use multiscale modeling techniques to reduce the dimensionality of the original problem \cite{efendiev2009multiscale, allaire2005multiscale, masud2006multiscale}. There are many modifications of multiscale methods where each may suit a specific task. The most famous is the multiscale finite element method (MsFEM) \cite{efendiev2009multiscale, hou1997multiscale, masud2006multiscale}. The MsFEM can give a significant error in domains with high-contrast properties. The Generalized Multiscale Finite Element Method (GMsFEM) is based on the spectral properties of the local problem and provides accurate approximation by defining multiple basis functions in local domains \cite{akkutlu2016multiscale, chung2017coupling, efendiev2013generalized, spiridonov2019generalized, spiridonov2020generalized}. Based on the finite volume method, a multiscale finite volume method (MsFVM) was developed \cite{jenny2005adaptive, sokolova2019multiscale}. 
The mixed finite element method(Mixed-FEM) has been developed to solve fluid flow problems in porous media \cite{chung2016mixed, auricchio2017mixed, boffi2013mixed, spiridonov2022mixed}.  
Multiscale methods with an oversampling strategy are used for problems with complex heterogeneity, such as channels or fractures. In the constraint energy minimizing generalized multiscale finite element method (CEM-GMsFEM), the multiscale basis functions are constructed in oversampled local domains and therefore take into account the influence of heterogeneity in neighboring local domains \cite{chung2018non, vasilyeva2019constrained, vasilyeva2019nonlocal, chung2018constraint, chung2018constraint1, cheung2020constraint}. For nonlinear problems, the online generalized multiscale finite method (Online GMsFEM) can be used where additional multiscale bases take into account changes in properties in nonlinear problems \cite{chung2015residual, chung2017online, spiridonov2023online, spiridonov2021online}.
In our previous work \cite{spiridonov2019mixed}, we used the mixed generalized finite element method for the Darcy-Fochheimer model in a heterogeneous domain. This paper extends the Darcy-Forchheimer model by adding fractures and time. As a multiscale method, we chose the non-local multicontinuum method (NLMC) \cite{vasilyeva2019upscaling, vasilyeva2019nonlocal, zhao2020analysis}. 

%
%

In this paper, we look at a non-local multicontinuum method (NLMC) for Darcy-Forchheimer flow in a fractured domain. The algorithm is divided into two parts: offline and online stages. In the offline stage, we construct multiscale basis functions. In this algorithm, we calculate bases in oversampled local domains with energy-minimizing constraints. At the online step, using the acquired bases, we solve the system on a coarse grid. It should be noted that the approximation relies on the finite volume method, and that fractures are represented by an embedded fracture model of lower dimensional fractures. Numerical results are presented for a two-dimensional fractured heterogeneous domain. The numerical experiment comprises examining the accuracy of the NLMC approach in relation to nonlinearity.


The Darcy-Forchheimer model is shown in a fractured heterogeneous domain in Section 2. A fine grid approximation using the finite volume method and embedded fracture model (EFM) is described in Section 3. The non-local multicontinuum technique (NLMC) algorithm for the Darcy-Forchheimer model is presented in Section 4. We describe multiscale basis functions in oversampled domains with constraints. In Section 5, we analyse the impact of nonlinearity on the method's accuracy and offer a numerical experiment for two test scenarios.


\section{Problem formulation}

The filtration process is described by the well-known Darcy law equation:
\begin{equation} \label{eq1}
\mu k^{-1} \boldsymbol{u}+\nabla p = 0.
\end{equation}
Let us supplement the basic equation of Darcy's law with a nonlinear term. The nonlinear Darcy-Forchheimer equation writing as follows:
\begin{equation} \label{eq2}
\mu k^{-1} \boldsymbol{u} + \rho \beta |\boldsymbol{u}|\boldsymbol{u} + \nabla p = 0,
\end{equation}
where $p$ is the pressure, $\boldsymbol{u}$ is the velocity, $k$ is the heterogeneous permeability, $\mu$ is the viscosity, $\rho$ is the density and $\beta$ is the Forchheimer coefficient. The impact of nonlinearity on the overall physical process is determined by the Forchheimer coefficient.

We can express the equation \eqref{eq2} as:
\begin{equation} \label{eq3}
\boldsymbol{u} = - \frac{k / \mu}{1+\rho \beta k |\boldsymbol{u}|/ \mu} \nabla p.
\end{equation}


In this study, we take a dynamic nonlinear filtration process into consideration and include a temporal derivative 
\begin{equation}\label{eq4}
\frac{\partial(\rho \phi)}{\partial t} + \nabla \cdot (\rho \boldsymbol{u}) = \rho f.
\end{equation}

Assuming $\rho$ as the constant, we set
\begin{equation} 
\frac{\partial \phi}{\partial t} = c_r \phi _0 \frac{\partial p}{\partial t}, \quad c_r=\frac{1}{\phi_0}\frac{\partial \phi}{\partial p}, \quad c=c_r\phi _0,
\end{equation} 
where $\phi$ is the porosity and $c_r$ is the porous media compressibility coefficient.

Therefore, the dynamic incompressible single-phase Darcy-Forchheimer flow equation can be written as follows:
\begin{equation} \label{eq5}
c\frac{\partial p}{\partial t} + \nabla \cdot  \boldsymbol{u} = f,
\end{equation}
where $f$ is a source term. 

In this paper, we consider the nonlinear filtration process in fractured media.  
We define the porous matrix domain as $\Omega \in \mathcal{R}^d$ and lower dimensional fractures as $\gamma \in \mathcal{R}^{d-1}$. In our implementation, we consider problems in the two-dimensional domain $d=2$.  
We  consider the problem in fractured media using the following model:
\begin{equation} \label{eq6}
\begin{split}
c_m \frac{\partial p_m}{\partial t} + \nabla \cdot  \boldsymbol{u}_m + r_{mf}(p_m,p_f)= f_m, \quad \boldsymbol{x} \in \Omega, \quad t > 0
\\
c_f \frac{\partial p_f}{\partial t} + \nabla \cdot  \boldsymbol{u}_f - r_{mf}(p_m,p_f)= f_f, \quad \boldsymbol{x} \in \gamma, \quad t > 0
\end{split}
\end{equation}
where $m$ and $f$ denotes subindices for matrix and fracture,  $r_{mf}, r_{fm}$ are the transfer terms between matrix-fracture and fracture-matrix. 

From  \eqref{eq3} for the velocity $\boldsymbol{u_i}$, we have 
\begin{equation} \label{eq7}
\boldsymbol{u_i} = - \frac{k_i / \mu}{1+\rho \beta_i k_i |\boldsymbol{u_i}|/ \mu} \nabla p_i, \quad i=m,f.
\end{equation}

Finally, the nonlinear Darcy-Forchheimer flow  in a fractured domain is defined by the following coupled system of equations
\begin{equation} \label{eq6a}
\begin{split}
c_m \frac{\partial p_m}{\partial t} - \nabla \cdot  \left( 
\frac{k_m / \mu}{1+\rho \beta_m k_m |\boldsymbol{u_m}|/ \mu} \nabla p_m
 \right) + r_{mf}(p_m,p_f)= f_m, \quad \boldsymbol{x} \in \Omega, \quad t > 0\\
c_f \frac{\partial p_f}{\partial t} - \nabla \cdot  \left( 
\frac{k_f / \mu}{1+\rho \beta_f k_f |\boldsymbol{u_f}|/ \mu} \nabla p_f
 \right)  - r_{mf}(p_m,p_f)= f_f, \quad \boldsymbol{x} \in \gamma, \quad t > 0.
\end{split}
\end{equation}

We supplement the equation \eqref{eq6a} with the following boundary conditions
\begin{equation} \label{eq8}
\boldsymbol{u}_m\cdot \boldsymbol{n} = 0, \quad \boldsymbol{u}_f\cdot \boldsymbol{n} = 0,\quad \boldsymbol{x}  \in \partial \Omega, \quad t>0, 
\end{equation}
and given initial condition $p_m = p_f = p_0$ for $t = 0$.

\section{Fine grid approximation}

Next, we consider the fine grid approximation of the problem \eqref{eq6a}, \eqref{eq8}. We use the finite volume method with an embedded fracture model. We define structured triangular fine grid $\mathcal{T}_h$, which does not conform to fractures.   We build a separate mesh for fractures and denote it as $\mathcal{E}_h$. 
For the time-dependent problem,  we note $n$ as a number of time layers, $\tau$ as the time step, and $T_{max}=n\tau$ as the final time. 
To make an approximation of the fine grid, we define the element of the fine grid $\varsigma _i$ and the element of fracture mesh $\iota _l$. We consider $N^m_f$ and $N^f_f$ as the number of elements in fine grids of matrix and fracture, respectively. From that, we can write a fine grid for porous matrix domain as $\mathcal{T}_h=\bigcup^{N^m_f}_{i=1} \varsigma_i$ and for fractures as $\mathcal{E}=\bigcup^{N^f_f}_{l=1} \iota_l$. 

We have the following approximation on the fine grid:
\begin{equation}\label{eq9}
\begin{split}
c_m \frac{p^{n+1}_{m,i}-p^n_{m,i}}{\tau}|\varsigma_i|+ \sum_{j} T_{ij}^n (p_{m,i}^{n+1} - p_{m,j}^{n+1}) + 
+ \sum_{l} \sigma^n_{il} (p_{m,i}^{n+1} - p_{f,l}^{n+1})=f_{m,i} |\varsigma_i|, \quad \forall i = 1, \ N^m_f, \\
c_f \frac{p^{n+1}_{f,l}-p^n_{f,l}}{\tau}|\iota _l|+ \sum_{\eta} W_{l\eta}^n (p_{f,l}^{n+1} - p_{f,\eta}^{n+1})-
-\sum_{i} \sigma^n_{il}(p_{m,i}^{n+1}-p_{f,l}^{n+1}) = f_{f,l}|\iota _l|, \quad \forall l=1, \ N^f_f,
\end{split}
\end{equation}
with
\[
T^n_{ij}=Z_{ij}\varrho^n_{ij}, \quad \mbox{where} \quad 
Z_{ij}=\frac{k_{m,ij} \ |E_{ij}|}{\mu \ d^m_{ij}}, \quad \varrho^n_{ij}=\frac{1}{1+\rho \beta_{m,ij} k_{m,ij} |u^n_{m,ij}|/\mu}, \nonumber
\]\[
W_{l\eta}^n = X_{l\eta} w^n_{l\eta} \quad \mbox{where} \quad 
X_{l\eta} = \frac{k_{f,l\eta}}{\mu \ d^f_{l\eta}}, \quad w^n_{l\eta}=\frac{1}{1+\rho \beta_{f,l\eta} k_{f,l\eta} |u^n_{f,l\eta}|/\mu}, \nonumber
\]
where  
\[
k_{m,ij}=\frac{2}{1/k_{m,i}+1/k_{m,j}}, \quad 
\beta_{m,ij}=\frac{2}{1/\beta_{m,i}+1/\beta_{m,j}}, 
\]\[
k_{f,l\eta}=\frac{2}{1/k_{f,l}+1/k_{f,\eta}}, \quad 
\beta_{f,l\eta}=\frac{2}{1/\beta_{f,l}+1/\beta_{f,\eta}},
\]and $|E_{ij}|$ is the lenght of facet between cells $\varsigma _i$ and $\varsigma _j$, $d_{ij}$ is the distance between midpoint of cells $\varsigma _i$ and $\varsigma _j$, 
and $d^f_{l\eta}$ is the distance between midpoint of cells $\iota _l$ and $\iota _\eta$.

For transfer terms, we use the following  approximation
\[
\sigma^n _{il} = Y_{il} \eth_{il}, \quad \mbox{where} \quad Y_{il}=\frac{k^*_{il} |\gamma_{il}|}{\mu \ \theta_{il}}, \quad \eth_{il}=\frac{1}{1+\rho \beta^*_{il}k^*_{il} |u^n_{mf.il}|/\mu,}
\]
where
\[
\quad k^*_{il} = \frac{2}{1/k_{m,i}+1/k_{f,l}}, \quad \beta^*_{il} = \frac{2}{1/\beta_{m,i}+1/\beta_{f,l}}, 
\]
and $C_{il}=\frac{|\gamma_{il}|}{\theta_{il}}$ is the connectivity index, $\theta_{il}$ is the distance midpoint of matrix cell $\varsigma _i$ and fracture cell $\iota _l$ ,  $|\gamma_{il}|$ is the length of the intersection of the fractures cell $\iota _l$ and matrix cell $\varsigma _i$.

Therefore, we can present the system of equations \eqref{eq9} in matrix form:
\begin{equation}\label{eq10}
M \frac{p^{n+1}-p^n}{\tau} + (A^n+Q^n)p^{n+1} = F, \quad \boldsymbol{x} \in \Omega,
\end{equation}
where
\begin{equation}
M=\begin{pmatrix}
M_m & 0 \\
0 & M_f
\end{pmatrix}, \quad A^n = \begin{pmatrix}
A^n_m & 0 \\
0 & A^n_f
\end{pmatrix}, \quad 
Q^n =\begin{pmatrix}
Q^n_{mf} & -Q^{n}_{mf} \\
-Q^n_{fm} & Q^n_{fm}
\end{pmatrix}, \quad F=\begin{pmatrix}
F_m \\
F_f
\end{pmatrix}, \nonumber
\end{equation}
and
\begin{equation}
\begin{split}
&M_m = \{m^m_{ij}\}, \quad 
m^m_{ij}= \begin{cases} c_m|\varsigma _i| & i=j, \\ 0 & i \neq j \\  \end{cases}, \quad M_f = \{m^f_{l\eta}\}, \quad m^f_{l\eta}= \begin{cases} c_f|\iota _l| & l=\eta, \\ 0 & l \neq \eta \\  \end{cases}, \\
&A^n_m = \{a^{m,n}_{ij}\}, \quad 
a^{m,n}_{ij}= \begin{cases} \sum_j T^n_{ij} & i=j, \\ -T^n_{ij} & i \neq j \\  \end{cases}, \quad A^n_f = \{a^{f,n}_{l\eta}\}, \quad a^{f,n}_{l\eta}= \begin{cases} \sum_\eta W^n_{l\eta} & l=\eta, \\ -W^n_{l\eta} & l \neq \eta \\  \end{cases}, \\
&Q^n_{fm} = \{q^{fm, n}_{il}\}, \quad 
q^{fm, n}_{ij}= \begin{cases} \sigma^n_{li} & \iota _l \subset \varsigma _i, \\ 0 & l \neq \eta \end{cases}, \quad 
Q^n_{mf} = (Q^n_{fm})^T, 
 \end{split}. 
\end{equation}

In this approximation, we have nonlinear terms $T^n_{ij} $, $W^n_{l\eta}$ and $ \sigma^n_{li}$ that can be factorized to the linear part $Z_{ij}$, $X_{l\eta}, Y_{il}$ and nonlinear part $\varrho^n_{ij}$, $w^n_{l\eta}$, $\eth_{il}$. We use the nonlinear velocity $\boldsymbol{u}_m$ and $\boldsymbol{u}_f$  from the previous time layer.


\section{Coarse grid approximation}

Next, we present a coarse grid approximation using the Non-Local Multi-Continuum method (NLMC). In this approach, we construct the multiscale basis functions in an oversampled local domain. To derive the multiscale basis function, we solve local problems with constraints. Each multiscale basis is calculated for one target continuum, determined by the imposed constraints. Multiscale basis functions use a constraint that makes an integral over the local domain (coarse cell) disappear in all continuums besides the target continuum and provide a meaning of the coarse grid solution. The basis functions were calculated for each fracture network rather than for each fracture individually. Such an approach of basis construction separates background and fractures and has spatial decay properties.

\begin{figure}[h!]
\begin{center}
\includegraphics[width=0.5\linewidth]{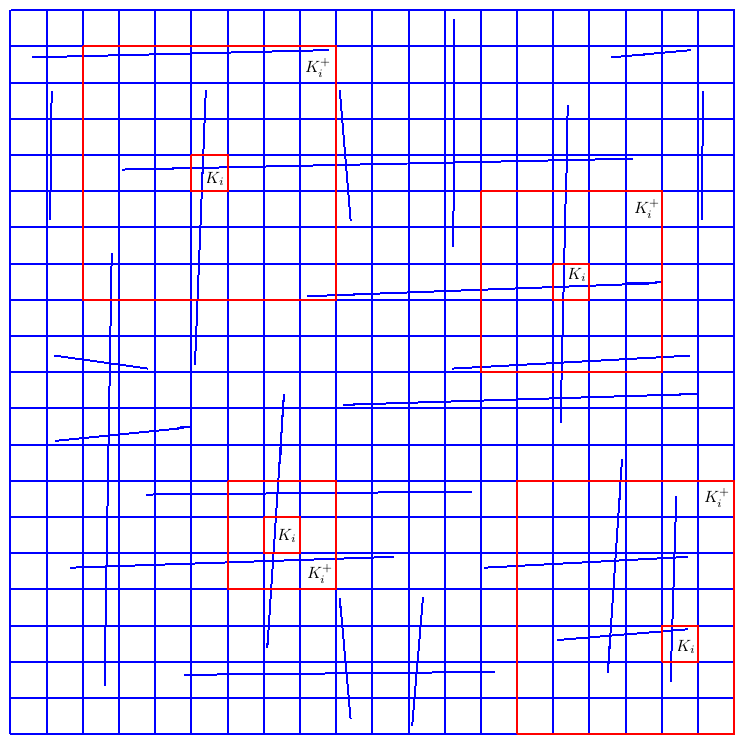}
\end{center}
\caption{Coarse grid with local domains $K^+_i$ and different numbers of oversampling layers.}
\label{localdom}
\end{figure}

Let $\mathcal{T}_H$ be the coarse grid divided into coarse cells $K_i$. We define an oversampled local domain $K^+_i$  by increasing $K_i$ by several coarse oversampling layers. We denote the number of oversampling layers as $S$. An example of local domains $K^+_i$ with a different number of oversampling layers is presented in Figure \ref{localdom}. The fracture network $\gamma^{(l)}$ connects into global fracture network $\gamma = \cup^L_{l=1}\gamma^{l}$, where $L$ is the total number of fracture networks in $K_i$. In each coarse cell $K_j$ we have $L_j$ amount of fracture network and $\gamma^{(l)}_j=K_j \cap \gamma^{(l)}$ is the fracture network laying inside $K_j$. In each coarse cell $K_i$, we construct $L_i + 1$ multiscale basis functions: $L_i$ functions for each fracture network $\gamma^{(l)}_i$ and one for background medium $K_i$.

The NLMC method is divided into two stages: offline and online. In the offline stage, we construct the multiscale basis functions. We solve a system on the coarse grid using multiscale basis functions in the online stage. 
To derive the multiscale basis functions, we solve the following problem in each oversampled local domain $K^+_i$:
\begin{equation} \label{eq11}
\begin{split}
 \nabla \cdot \left( \frac{k_m}{\mu} \nabla \psi^i_m   \right)+ \tilde{r}_{mf}(\psi^i_m-\psi^i_f)= 0, \\
 \nabla \cdot \left( \frac{k_f}{\mu} \nabla \psi^i_f  \right) -  \tilde{r}_{mf}(\psi^i_m-\psi^i_f)= 0, 
\end{split}
\end{equation}
with 
\[
 \psi^i_m = 0, \quad \psi^i_f = 0, \quad \boldsymbol{x} \in \partial K^+_i. 
 \]
Here we use a linear part of transfer term $r_{mf}$, which will be explained further below.

We apply the following constraints for the problem \eqref{eq11} in each $K_j \subset K^+_i $:

\begin{itemize}
\item for background:
\begin{equation} \label{eq12}
\frac{1}{|K_j|} \int_{K_j} \psi^{i, 0}_{m} dx = \delta_{i,j} , \quad 
\frac{1}{|\gamma^{(m)}_j|}\int _{\gamma^{(m)}_j} \psi^{i, 0}_{f} ds = 0, \quad m = \overline{1, L_j},
\end{equation}
\item for fractures:
\begin{equation} \label{eq13}
\frac{1}{|K_j|} \int_{K_j} \psi^{i,l}_{m} dx = 0 , \quad 
\frac{1}{|\gamma^{(m)}_j|} \int _{\gamma^{(m)}_j} \psi^{i,l}_{f} ds = \delta_{i,j} \delta_{m,l}, \quad m,l = \overline{1, L_j},
\end{equation}
\end{itemize}
with $\psi^{i,l}=(\psi^{i,0}_{m},\psi^{i,0}_{f})$ and $\psi^{i,l}=(\psi^{i,l}_{m},\psi^{i,l}_{f})$.

The constraints provide a function for background medium that has a mean value of one in $K_i$ and a mean value of zero elsewhere in $K^+_i$. Additionally, all fractures inside the oversampled local domain $K^+_i$ will have mean values of one for the multiscale basis function. We have a basis function for the fracture network in the background continuum that has a mean value of zero for each coarse cell in $K^+_i$. Moreover, for the target fracture network $\gamma ^{(l)}_i$, we have basis functions with mean values one and zero for all fracture networks inside $K^+_i$.

We should mention that we compute multiscale basis functions for the linear parts of $T_{ij}$ and $W_{l\eta}$ from \eqref{eq9}. The multiscale basis functions were computed only once and did not vary over time. Therefore, the approximation of the problem \eqref{eq11} takes the following form:
\begin{equation} \label{eq14}
\begin{split}
\sum_{j} Z_{kj} (\psi^{i,l}_{m,k} - \psi^{i,l}_{m,j}) + \sum_{n} Y _{kn} (\psi^{i,l}_{m,k} - \psi^{i,l}_{f,n})= 0, \quad \forall k = \overline{1, N^m_{f, K_i^+}} \\
\sum_{\eta} X_{l\eta} (\psi^{i,l}_{f,n} - \psi^{i,l}_{f,\eta})-\sum_{k} Y_{il}(\psi^{i,l}_{m,k}-\psi^{i,l}_{f,n}) = 0, \quad \forall n=\overline{1,{N^f_{f, K_i^+}}}.
\end{split}
\end{equation}

\begin{figure}[h!]
\begin{center}
\includegraphics[width=0.9\linewidth]{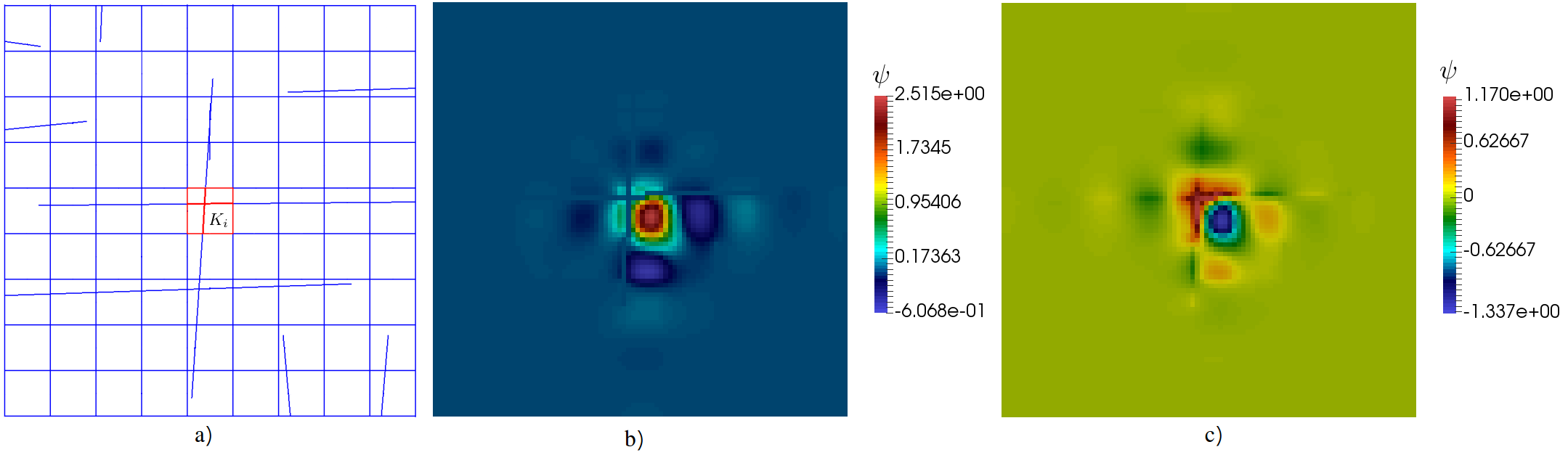}
\end{center}
\caption{a) Oversampled local domain $K^+_i$. b) Multiscale basis function for matrix c) Multiscale basis function for fracture network.}
\label{fig2}
\end{figure}

In approximation \eqref{eq14}, all notations are taken from \eqref{eq9}, except $N^m_{f, K_i^+}$ and $N^f_{f, K_i^+}$, which represent the number of fine grids elements in oversampled local domain $K^+_i$. 
The approximation of the local problem for the development of multiscale basis functions can be represented in the matrix form: 
\begin{equation} \label{eq15}
\begin{pmatrix}
A^i_m + Q^i_{mf} & -Q^i_{mf} & (B^i_m)^T & 0 \\
-Q^i_{fm} & A^i_f + Q^i_{fm} & 0 & (B^i_f)^T \\
B^i_m & 0 & 0 & 0 \\
0 & B^i_f & 0 & 0
\end{pmatrix} \begin{pmatrix}
\psi^{i,l} _m \\
\psi^{i,l} _f \\
\varphi^{i,l} _m \\
\varphi^{i,l} _f
\end{pmatrix} = \begin{pmatrix}
0 \\
0 \\
F^{i,l}_m \\
F^{i,l}_f
\end{pmatrix},  
\end{equation}
where $\varphi^i _m$ and $\varphi^i _f$ are Lagrange multipliers, that comes from constraints.  
For matrix continuum (first basis) we set  $F^{i,0}_m = \{f^{i,0}_{m,j}\}$, $ f^{i,0}_{m,j} = \delta_{i,j}$ and $F^{i,0}_f = \boldsymbol{0}$. 
For $m$th fracture we set $F^{i,l}_m = \boldsymbol{0}$ and $F^{i,l}_f = \{f^{i,l}_{f,j}\}$, $ f^{i,l}_{f,j} = \delta_{i,j}\delta _{m,l}$. Figure \ref{fig2} demonstrates an illustration of multiscale basis functions that was computed in $K^+_i$ using four oversampling layers.

Next, we can define the multiscale space using obtained multiscale basis functions:
\begin{equation}
V_{ms} = \mbox{span}\{(\psi^{i,l}_m, \psi^{i,l}_f), \ i=\overline{1, N_c}, \ l = \overline{0,L_i}\}.
\end{equation}

To solve the problem on the coarse grid, we define a projection matrix on multiscale space:
\begin{equation}
R  =\begin{pmatrix}
R_{mm} & R_{mf} \\
R_{fm} & R_{ff}
\end{pmatrix}, 
\end{equation}
\[
\begin{split}
R_{mm} & = [ \psi^{0,0}_m, \psi^{1,0}_m, ... , \psi^{N_c,0}_m]^T,   \quad
R_{ff}  = [ \psi^{0,1}_f, ..., \psi^{0,L_0}_f, \psi^{1,1}_f, ..., \psi^{1,L_1}_f, ..., \psi^{N_c,1}_f, ..., \psi^{N_c,L_{N_c}}_f]^T, \\
R_{fm}  &= [ \psi^{0,1}_m, ..., \psi^{0,L_0}_m, \psi^{1,1}_m, ..., \psi^{1,L_1}_m, ..., \psi^{N_c,1}_m, ..., \psi^{N_c,L_{N_c}}_m]^T, \quad 
R_{mf} = [ \psi^{0,0}_f, \psi^{1,0}_f, ... , \psi^{N_c,0}_f]^T.
\end{split} 
\]

Next, we  write approximation on the coarse grid for $\overline{p}^{n+1} = (\overline{p}^{n+1}_m, \overline{p}^{n+1}_f)^T$
\begin{equation} \label{eq16}
\overline{M} \frac{\overline{p}^{n+1}-\overline{p}^n}{\tau}
+ 
(\overline{A}^n +  \overline{Q}^n)\overline{p}^{n+1} = \overline{F},
\end{equation}
where 
\begin{equation}
\overline{A}^n = RA^nR^T, \quad 
\overline{M} = RMR^T, \quad 
\overline{Q}^n = RQ^nR^T, \quad \overline{F}=RF.
\end{equation}

The resulting functions $\overline{p}_m^{n+1}$ and $\overline{p}_f^{n+1}$ store the average values over the coarse grid element. In addition, we can reconstruct a solution on the fine grid by $p_{ms}^{n+1}=R^T\overline{p}^{n+1}$.

Finally, we can summarize the algorithm of the NLMC method as follows:
\begin{enumerate}
\item Define the coarse grid $\mathcal{T}_H$ and  oversampled local domains $K^+_i$ for a given number of oversampling layers. 
\item Solve local problems \eqref{eq15} in each $K^+_i$ to get the multiscale basis functions;
\item Construct the projection matrix $R$ by combining multiscale basis functions;
\item Project and solve the nonlinear system on the coarse grid \eqref{eq16}.
\end{enumerate}

\section{Numerical results}

\begin{figure}[h!]
\begin{center}
\includegraphics[width=0.9\linewidth]{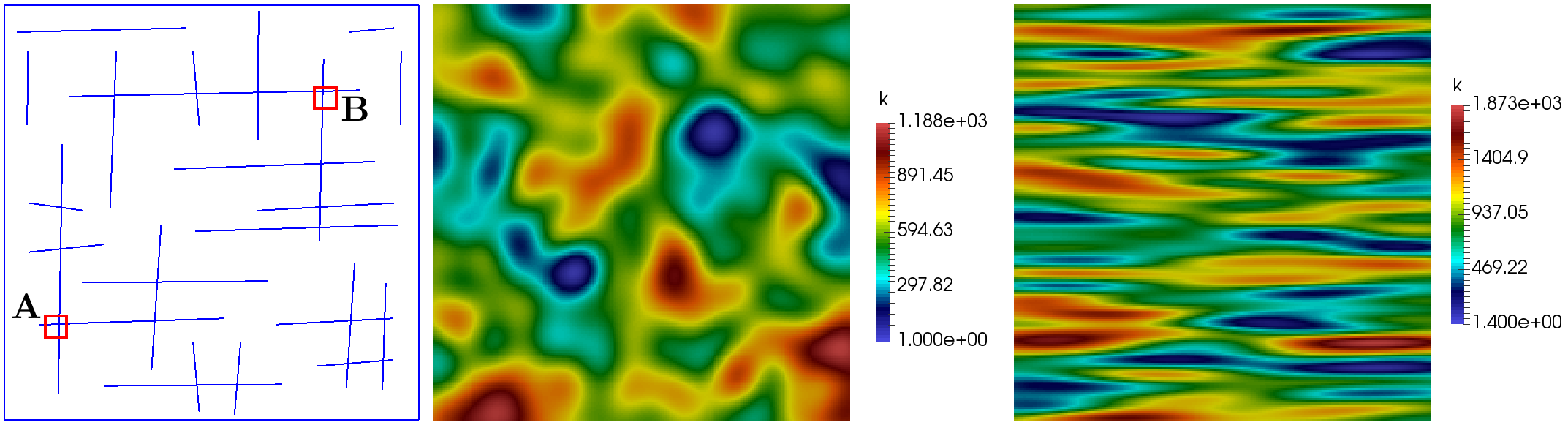}
\end{center}
\caption{Heterogeneity and fracture distribution. Left: computational domain $\Omega$ with fractures $\gamma$. Middle: Coefficient $k_m$ in \textit{Test 1}. Right: Coefficient $k_m$ in \textit{Test 2}.}
\label{domkle}
\end{figure}

In this section, we use NLMC to provide numerical results for the Darcy-Forchheimer problem. We consider the flow problem in the two-dimensional fractured heterogeneous domain. Figure \ref{domkle} presents the computation domain with the fracture location. We present numerical results for two test cases with different heterogeneous coefficients $k_m$ and  denote test problems as \textit{Test 1} and \textit{Test 2}. We show the heterogeneous coefficient $k_m$ in Figure \ref{domkle}. We chose a highly heterogeneous domain for \textit{Test 2} with a significant contrast in $k_m$ values inside one coarse cell $K_i$. For fractures, we set $k_f=10^9$ to test both cases. 

We consider numerical experiments for two coarse grids: $20 \times 20$ and $40 \times 40$. Simulations are performed for $T_{max}=12.5 \cdot 10^5 $ with 100 time layers. The large number of time layers is used to resolve the nonlinearity using an approximation from the previous time layer. To investigate the accuracy of the NLMC algorithm, we compare the multiscale solutions with the reference solutions. We take a fine-grid solution as a reference solution. We use  $200 \times 200$ uniform mesh with square cells for fine grid solution in a porous medium. For the fracture network, we use a one-dimensional mesh with 1730 elements. The fine grid system's size equals $41733 \times 41733$. 
We compare solutions by using $L_2$ norm errors:
\begin{equation} \label{eq17}
||e||^{L_2}_{p}=\sqrt{\frac{\int_{\Omega}(p_f-p_{ms})^2 dx}{\int_{\Omega}p^2_f dx}},  \quad ||\overline{e}||^{L_2}_{p}=\sqrt{\frac{\int_{\Omega}(\overline{p}_f-\overline{p}_{ms})^2 dx}{\int_{\Omega}\overline{p}^2_f dx}},
\end{equation}     
where $p_{ms}$ is the multiscale solution, $p_f$ is the solution on the fine grid, $\overline{p}_{ms}$, $\overline{p}_{f}$ are the average values over the coarse grid element for multiscale and fine grid solutions, respectively.

\begin{figure}[h!]
\begin{center}
\includegraphics[width=0.9\linewidth]{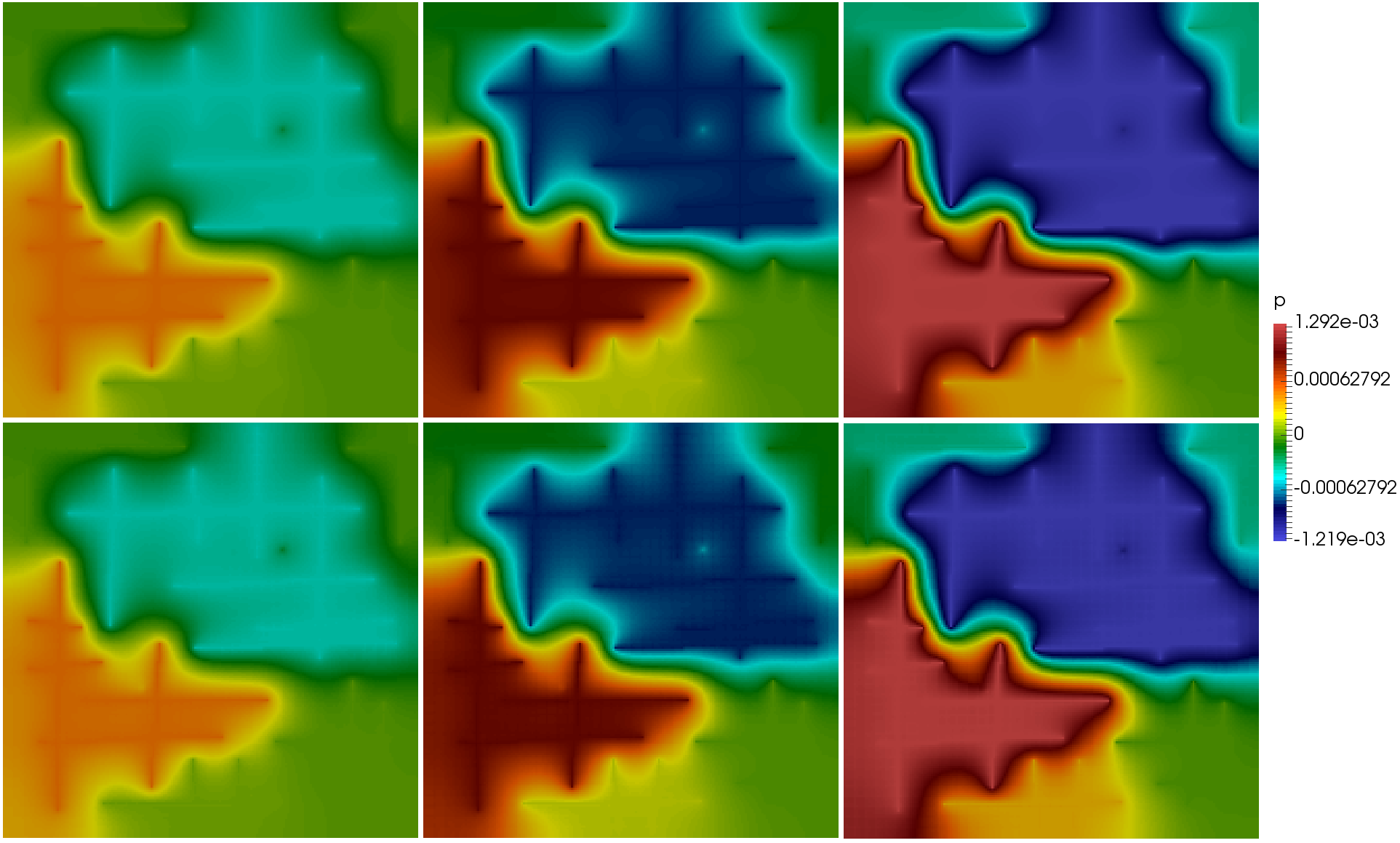}
\end{center}
\caption{Numerical results for $C = 10^4$ at the 30th, 60th, and last time layer. First row: fine grid solution. Second row: NLMC-solution on $40 \times 40$ coarse grid using four oversampling layers. \textit{Test 1}.}
\label{results1}
\end{figure}

\begin{figure}[h!]
\begin{center}
\includegraphics[width=0.9\linewidth]{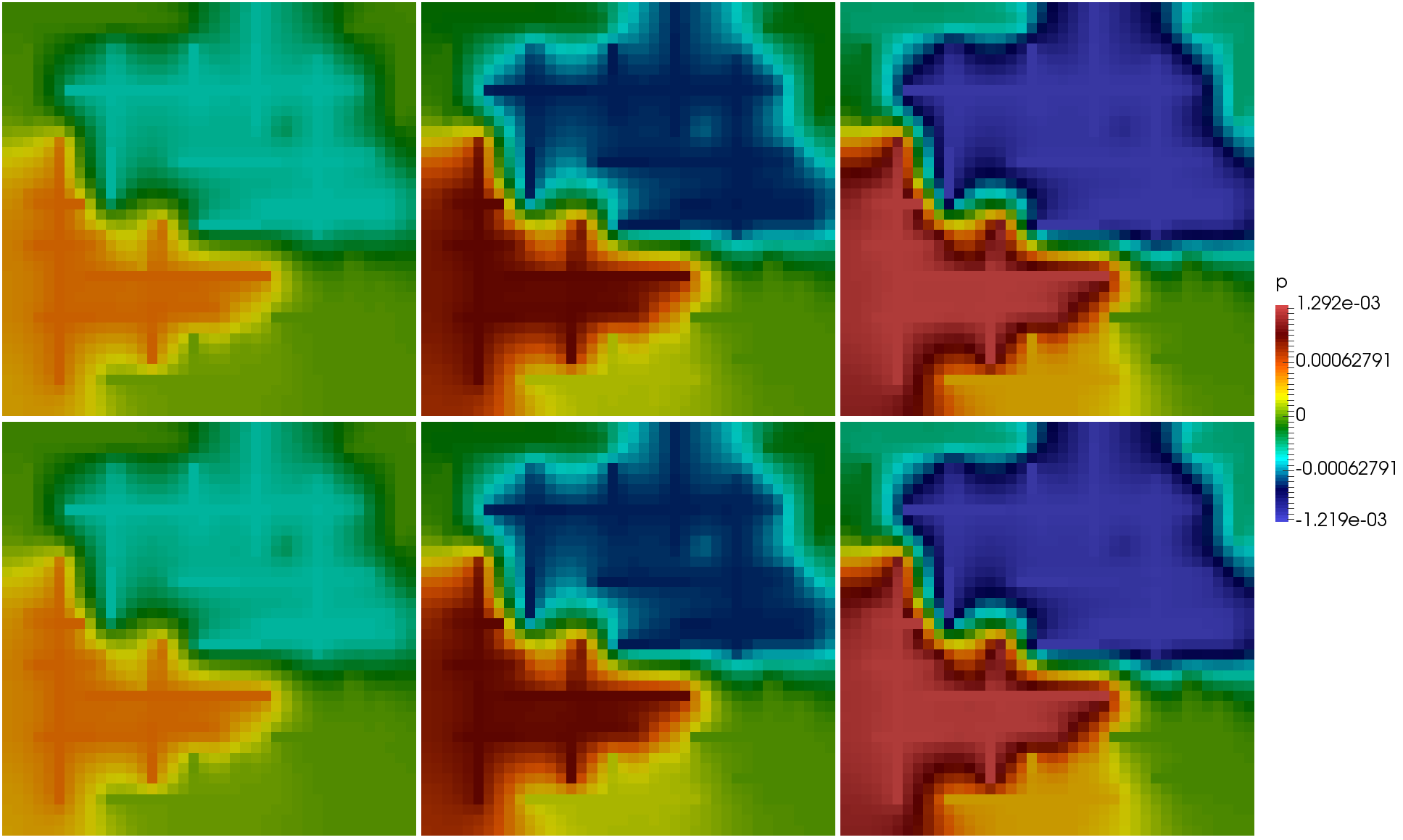}
\end{center}
\caption{Numerical results for $C = 10^4$ at the 30th, 60th, and last time layer. First row: Coarse-grid average of the fine grid solution on $40 \times 40$ coarse grid. Second row:  NLMC-solution on $40 \times 40$ coarse grid using four oversampling layers. \textit{Test 1}.}
\label{results1mean}
\end{figure}

We consider results with different values of $\beta_i$ to investigate the effect of nonlinearity. In  implementation, we set $\beta_m=Ck_m^{-1}$ and  $\beta_f=Ck_f^{-1}$. We change the value of $C$ to control the effect of nonlinearity. We consider results for $C=10^4$, $10^3$, $10^2$, $10$, $0$. We set $\mu=8$, $\rho=1.0$, $c_m=c_f=1.0$, $r_{mf}=k_m$, $f_m=0$. The source term for fractures $f_f$ contains two wells: injection and production. We put injection well in all fracture cells inside coarse cell \textbf{A} with $f_A = 10^-3$ and production well in all fracture cells inside coarse cell \textbf{B} with $f_B = -10^-3$. The location of injection and production wells are presented in Figure \ref{domkle}.   

\begin{figure}[h!]
\begin{center}
\includegraphics[width=0.9\linewidth]{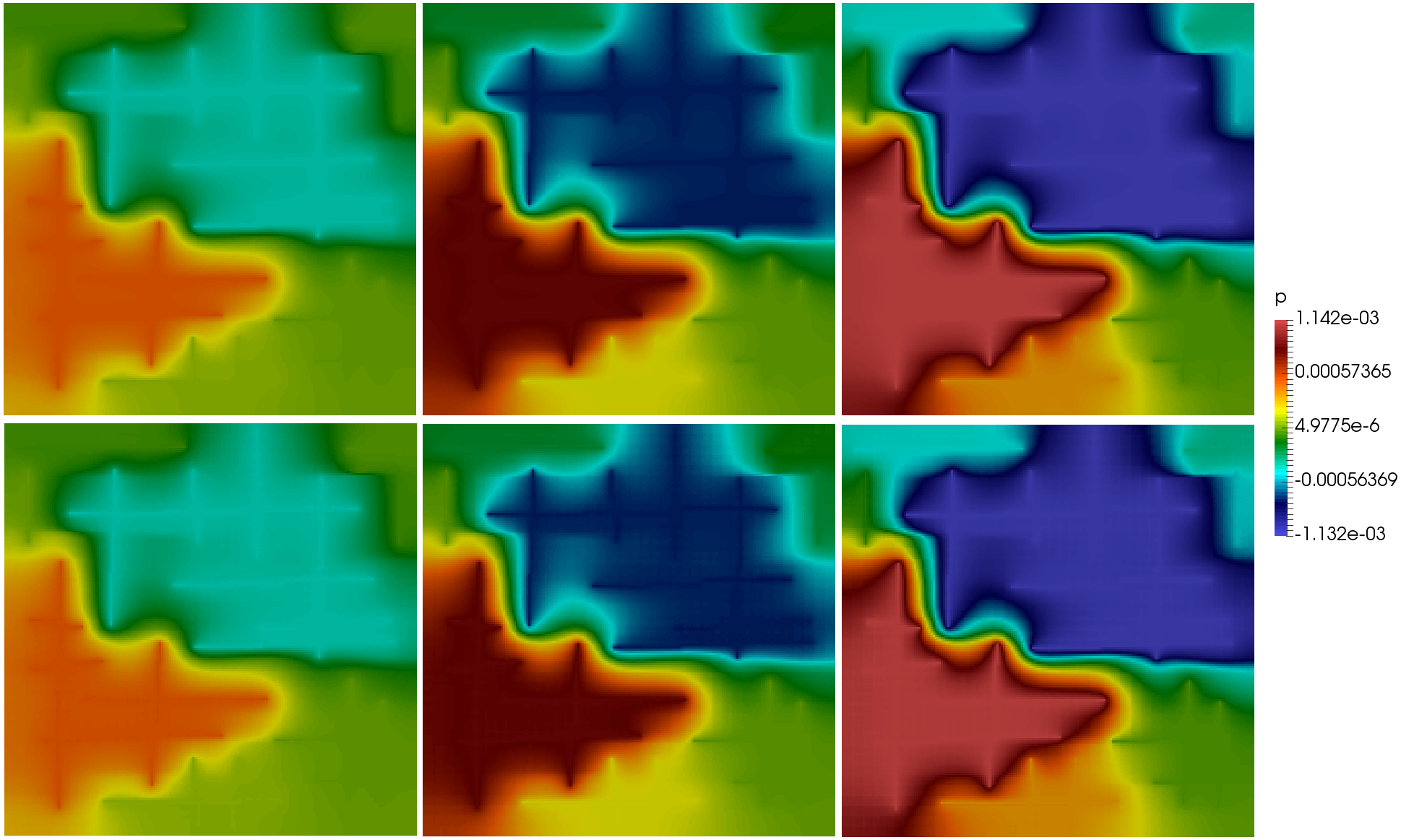}
\end{center}
\caption{Numerical results for $C = 10^4$ at the 30th, 60th, and last time layers. First row: fine grid solution. Second row: NLMC-solution on $40 \times 40$ coarse grid using four oversampling layers. \textit{Test 2}.}
\label{results2}
\end{figure}

\begin{figure}[h!]
\begin{center}
\includegraphics[width=0.9\linewidth]{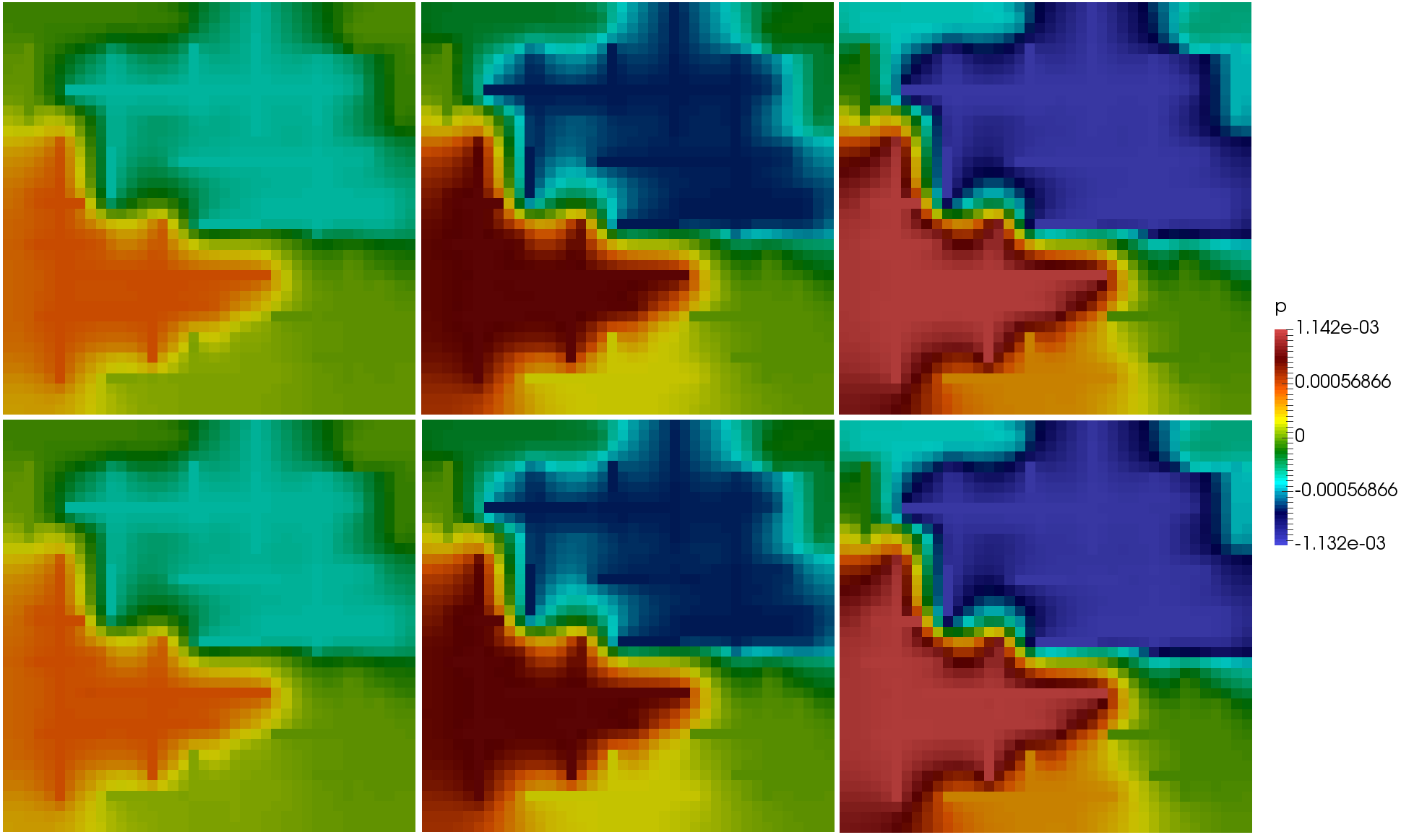}
\end{center}
\caption{Numerical results for $C = 10^4$ at the 30th, 60th, and last time layers. First row: Coarse-grid average of the fine grid solution on $40 \times 40$ coarse grid. Second row:  NLMC-solution on $40 \times 40$ coarse grid using four oversampling layers. \textit{Test 2}.}
\label{results2mean}
\end{figure}

The numerical results for \textit{Test 1} is presented in Figure \ref{results1} for $C=10^4$. Figure \ref{results1mean} shows an  solution on $40\times 40$ coarse grid for \textit{Test 1} with  $C=10^4$. The first row shows a solution on a fine grid, while the second row shows a solution using NLMC with four oversampling layers in basis construction. We present solutions for the 30th, 60th, and final time layers to demonstrate time distribution.  In the figures, the NLMC solution and the fine grid solution are looks similar. 

For \textit{Test 2}, we present results in Figure \ref{results2}. Figure \ref{results2mean} shows  solution on $40\times 40$ coarse grid for \textit{Test 2} with  $C=10^4$. The results are presented in the same order as \textit{Test 1} for $C=10^4$ and four oversampling layers in NLMC. We observe the same behavior as  \textit{Test 1}. The results are presented for the case with a larger influence of nonlinearity, $C=10^4$. 
 
\begin{table}[h!]
\begin{center}
\begin{tabular}{|c|c|c|c|c|}
\hline
  &  \multicolumn{2}{|c|}{$20 \times 20$}  &  \multicolumn{2}{|c|}{$40 \times 40$}   \\
$S$ &   \multicolumn{2}{|c|}{Coarse grid} &  \multicolumn{2}{|c|}{Coarse grid} \\
\cline{2-5}
  & $e^{L_2}_p$, (\%) & $\overline{e}^{L_2}_p$, (\%) & $e^{L_2}_p$, (\%) & $\overline{e}^{L_2}_p$, (\%) \\
\hline
3  & 6.583 & 6.337 & 5.576 & 5.364 \\
4  & 0.826 & 0.345 & 0.693 & 0.403 \\
5  & 0.179 & 0.017 & 0.144 & 0.029 \\
6  & 0.039 & 0.001 & 0.035 & 0.001\\
7  & 0.014 & 0.001 & 0.009 & 0.001 \\
\hline
\end{tabular}
\begin{tabular}{|c|c|c|c|c|}
\hline
   &  \multicolumn{2}{|c|}{$20 \times 20$}  &  \multicolumn{2}{|c|}{$40 \times 40$}   \\
$S$ &   \multicolumn{2}{|c|}{Coarse grid} &  \multicolumn{2}{|c|}{Coarse grid} \\
\cline{2-5}
  & $e^{L_2}_p$, (\%) & $\overline{e}^{L_2}_p$, (\%) & $e^{L_2}_p$, (\%) & $\overline{e}^{L_2}_p$, (\%) \\
\hline
3  & 8.165 & 7.922 & 6.674 & 6.643 \\
4  & 0.944 & 0.487 & 0.751 & 0.503 \\
5  & 0.198 & 0.026 & 0.144 & 0.036 \\
6  & 0.046 & 0.001 & 0.035 & 0.001 \\
7  & 0.013 & 0.001 & 0.009 & 0.001 \\
\hline
\end{tabular}
\end{center}
\caption{Numerical results for $C = 0$. Left: \textit{Test 1}. Right: \textit{Test 2}.}
\label{table1}
\end{table} 

 
\begin{table}[h!]
\begin{center}
\begin{tabular}{|c|c|c|c|c|}
\hline
  &  \multicolumn{2}{|c|}{$20 \times 20$}  &  \multicolumn{2}{|c|}{$40 \times 40$}   \\
$S$ &   \multicolumn{2}{|c|}{Coarse grid} &  \multicolumn{2}{|c|}{Coarse grid} \\
\cline{2-5}
  & $e^{L_2}_p$, (\%) & $\overline{e}^{L_2}_p$, (\%) & $e^{L_2}_p$, (\%) & $\overline{e}^{L_2}_p$, (\%) \\
\hline
3  & 6.583 & 6.338 & 5.577 & 5.365 \\
4  & 0.826 & 0.345 & 0.693 & 0.403 \\
5  & 0.179 & 0.017 & 0.144 & 0.029\\
6  & 0.039 & 0.001 & 0.035 & 0.001 \\
7  & 0.014 & 0.001 & 0.009 & 0.001\\
\hline
\end{tabular}
\begin{tabular}{|c|c|c|c|c|}
\hline
   &  \multicolumn{2}{|c|}{$20 \times 20$}  &  \multicolumn{2}{|c|}{$40 \times 40$}   \\
$S$ &   \multicolumn{2}{|c|}{Coarse grid} &  \multicolumn{2}{|c|}{Coarse grid} \\
\cline{2-5}
  & $e^{L_2}_p$, (\%) & $\overline{e}^{L_2}_p$, (\%) & $e^{L_2}_p$, (\%) & $\overline{e}^{L_2}_p$, (\%) \\
\hline
3  & 8.165 & 7.923 & 6.674 & 6.643 \\
4  & 0.944 & 0.487 & 0.749 & 0.503 \\
5  & 0.198 & 0.026 & 0.144 & 0.036 \\
6  & 0.046 & 0.001 & 0.035 & 0.001\\
7  & 0.013 & 0.001 & 0.009 & 0.001 \\
\hline
\end{tabular}
\end{center}
\caption{Numerical results for $C = 10$. Left: \textit{Test 1}. Right: \textit{Test 2}.}
\label{table2}
\end{table}


\begin{table}[h!]
\begin{center}
\begin{tabular}{|c|c|c|c|c|}
\hline
  &  \multicolumn{2}{|c|}{$20 \times 20$}  &  \multicolumn{2}{|c|}{$40 \times 40$}   \\
$S$ &   \multicolumn{2}{|c|}{Coarse grid} &  \multicolumn{2}{|c|}{Coarse grid} \\
\cline{2-5}
  & $e^{L_2}_p$, (\%) & $\overline{e}^{L_2}_p$, (\%) & $e^{L_2}_p$, (\%) & $\overline{e}^{L_2}_p$, (\%) \\
\hline
3  & 6.583 & 6.339 & 5.577 & 5.365 \\
4  & 0.826 & 0.345 & 0.693 & 0.403 \\
5  & 0.181 & 0.017 & 0.144 & 0.031 \\
6  & 0.039 & 0.001 & 0.035 & 0.001 \\
7  & 0.013 & 0.001 & 0.009 & 0.001 \\
\hline
\end{tabular}
\begin{tabular}{|c|c|c|c|c|}
\hline
   &  \multicolumn{2}{|c|}{$20 \times 20$}  &  \multicolumn{2}{|c|}{$40 \times 40$}   \\
$S$ &   \multicolumn{2}{|c|}{Coarse grid} &  \multicolumn{2}{|c|}{Coarse grid} \\
\cline{2-5}
  & $e^{L_2}_p$, (\%) & $\overline{e}^{L_2}_p$, (\%) & $e^{L_2}_p$, (\%) & $\overline{e}^{L_2}_p$, (\%) \\
\hline
3  & 8.163 & 7.924 & 6.675 & 6.645 \\
4  & 0.944 & 0.487 & 0.751 & 0.503 \\
5  & 0.198 & 0.026 & 0.144 & 0.036 \\
6  & 0.046 & 0.001 & 0.035 & 0.001 \\
7  & 0.014 & 0.001 & 0.009 & 0.001 \\
\hline
\end{tabular}
\end{center}
\caption{Numerical results for $C = 10^2$. Left: \textit{Test 1}. Right: \textit{Test 2}.}
\label{table3}
\end{table}


\begin{table}[h!]
\begin{center}
\begin{tabular}{|c|c|c|c|c|}
\hline
  &  \multicolumn{2}{|c|}{$20 \times 20$}  &  \multicolumn{2}{|c|}{$40 \times 40$}   \\
$S$ &   \multicolumn{2}{|c|}{Coarse grid} &  \multicolumn{2}{|c|}{Coarse grid} \\
\cline{2-5}
  & $e^{L_2}_p$, (\%) & $\overline{e}^{L_2}_p$, (\%) & $e^{L_2}_p$, (\%) & $\overline{e}^{L_2}_p$, (\%) \\
\hline
3  & 6.584 & 6.349 & 5.582 & 5.372 \\
4  & 0.826 & 0.346 & 0.693 & 0.403 \\
5  & 0.179 & 0.017 & 0.144 & 0.029 \\
6  & 0.039 & 0.001 & 0.035 & 0.001 \\
7  & 0.015 & 0.001 & 0.009 & 0.001 \\
\hline
\end{tabular}
\begin{tabular}{|c|c|c|c|c|}
\hline
   &  \multicolumn{2}{|c|}{$20 \times 20$}  &  \multicolumn{2}{|c|}{$40 \times 40$}   \\
$S$ &   \multicolumn{2}{|c|}{Coarse grid} &  \multicolumn{2}{|c|}{Coarse grid} \\
\cline{2-5}
  & $e^{L_2}_p$, (\%) & $\overline{e}^{L_2}_p$, (\%) & $e^{L_2}_p$, (\%) & $\overline{e}^{L_2}_p$, (\%) \\
\hline
3  & 8.145 & 7.935 & 6.684 & 6.659 \\
4  & 0.944 & 0.489 & 0.751 & 0.503 \\
5  & 0.199 & 0.026 & 0.144 & 0.036 \\
6  & 0.047 & 0.001 & 0.035 & 0.001 \\
7  & 0.015 & 0.001 & 0.009 & 0.001 \\
\hline
\end{tabular}
\end{center}
\caption{Numerical results for $C = 10^3$. Left: \textit{Test 1}. Right: \textit{Test 2}.}
\label{table4}
\end{table}

\begin{table}[h!]
\begin{center}
\begin{tabular}{|c|c|c|c|c|}
\hline
  &  \multicolumn{2}{|c|}{$20 \times 20$}  &  \multicolumn{2}{|c|}{$40 \times 40$}   \\
$S$ &   \multicolumn{2}{|c|}{Coarse grid} &  \multicolumn{2}{|c|}{Coarse grid} \\
\cline{2-5}
  & $e^{L_2}_p$, (\%) & $\overline{e}^{L_2}_p$, (\%) & $e^{L_2}_p$, (\%) & $\overline{e}^{L_2}_p$, (\%) \\
\hline
3  & 6.649 & 6.489 & 5.647 & 5.458 \\
4  & 0.837 & 0.365 & 0.693 & 0.405 \\
5  & 0.193 & 0.017 & 0.145 & 0.032 \\
6  & 0.075 & 0.004 & 0.038 & 0.001 \\
7  & 0.065 & 0.004 & 0.018 & 0.001 \\
\hline
\end{tabular}
\begin{tabular}{|c|c|c|c|c|}
\hline
   &  \multicolumn{2}{|c|}{$20 \times 20$}  &  \multicolumn{2}{|c|}{$40 \times 40$}   \\
$S$ &   \multicolumn{2}{|c|}{Coarse grid} &  \multicolumn{2}{|c|}{Coarse grid} \\
\cline{2-5}
  & $e^{L_2}_p$, (\%) & $\overline{e}^{L_2}_p$, (\%) & $e^{L_2}_p$, (\%) & $\overline{e}^{L_2}_p$, (\%) \\
\hline
3  & 8.188 & 8.111 & 6.822 & 6.860 \\
4  & 0.961 & 0.527 & 0.752 & 0.513 \\
5  & 0.219 & 0.024 & 0.145 & 0.037 \\
6  & 0.094 & 0.008 & 0.039 & 0.001 \\
7  & 0.084 & 0.008 & 0.022 & 0.001 \\
\hline
\end{tabular}
\end{center}
\caption{Numerical results for $C = 10^4$. Left: \textit{Test 1}. Right: \textit{Test 2}.}
\label{table5}
\end{table}

We present a relative $L_2$ error in Tables \ref{table1}-\ref{table5}. We give results with varying numbers of oversampling layers, $S=3,4,5,6$ and $7$. In the left table, the errors for \textit{Test 1} are presented. The results for \textit{Test 2} are presented in the right table.  
In the tables, we demonstrate two errors: $e^{L_2}_2$ is the error on the fine grid and $\overline{e}^{L_2}_2$ is the error determined on the coarse grid. 
We consider accuracy for two coarse grids: $20 \times 20$ and $40 \times 40$. The size of the coarse system is equal to $562 \times 562$ for $20 \times 20$ coarse grid and $1930 \times 1930$ for $40 \times 40$ coarse grid. We see that the size of the coarse system is smaller than the fine grid system discussed above. 
Moreover, we can see that employing four oversampling layers in all computations is enough to produce an accurate solution. The solution diverges when we use 1 or 2 oversampling layers, and we did not include these results in the tables. The multiscale basis functions with three layers provide a solution with low accuracy. We have a small difference in error between \textit{Test 1} and \textit{Test 2}
We see that the difference in accuracy is greater when we employ a small number of oversampling layers. However, it becomes very small when we increase the number of oversampling layers. This points out that the NLMC method accuracy is practically independent of the type of heterogeneity. 
Also, we investigatr the solution's accuracy to the equation's nonlinear part. We have a very small increase in error when we are raising the value of $C$, including linear case with $C=0$.  We can see a difference in error between $C=10^4$ and $C=10^3$, but the error is nearly similar in other values of $C$.It demonstrates that the NLMC only responds to large changes in nonlinearity, but it produces very good results even at bigger nonlinearity effects. We also analyze the method's accuracy to the coarse mesh size. The accuracy is better in coarse grid $40 \times 40$ than in coarse grid $20 \times 20$. 

\begin{figure}[h!]
\begin{center}
\includegraphics[width=0.45\linewidth]{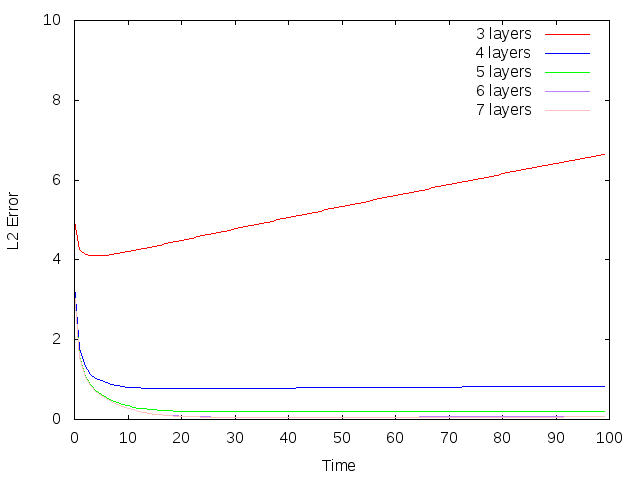}
\includegraphics[width=0.45\linewidth]{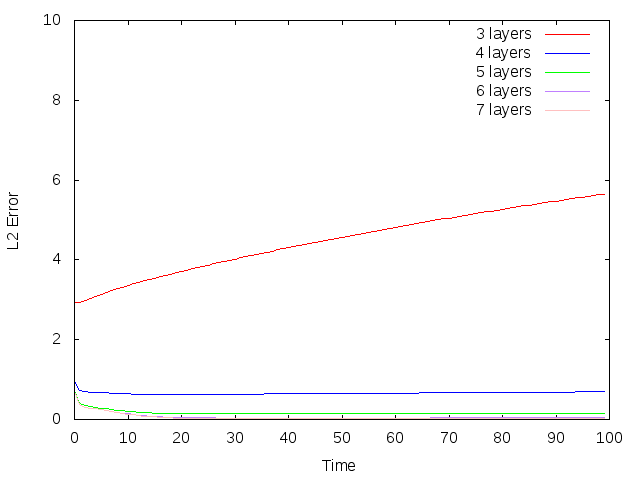}
\end{center}
\caption{\textit{Test 1}. Numerical results  for $C = 10^4$. Relative $L_2$ error distribution in time. Left: $20 \times 20$ coarse grid. Right: $40 \times 40$ coarse grid. }
\label{errorsC4}
\end{figure}

\begin{figure}[h!]
\begin{center}
\includegraphics[width=0.45\linewidth]{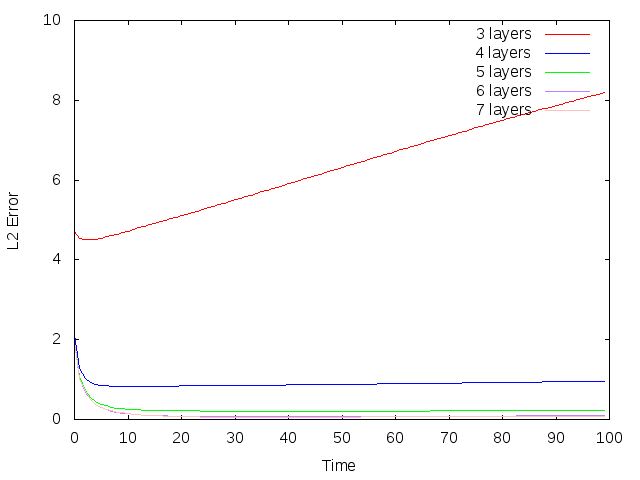}
\includegraphics[width=0.45\linewidth]{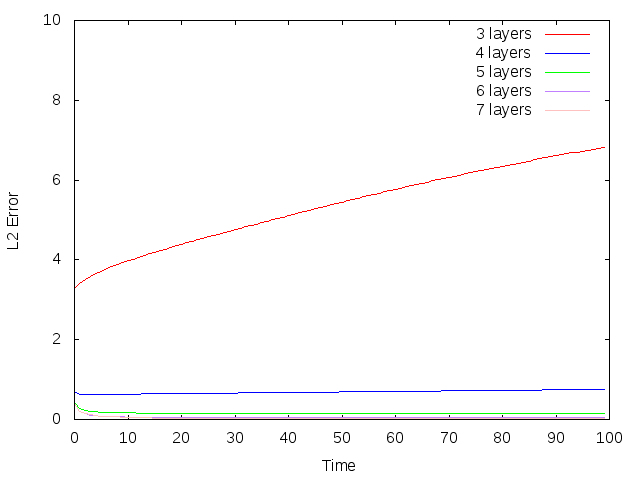}
\end{center}
\caption{\textit{Test 2}. Numerical results  for $C = 10^4$. Relative $L_2$ error distribution in time. Left: $20 \times 20$ coarse grid. Right: $40 \times 40$ coarse grid. }
\label{errorsC4_2}
\end{figure}

To show the distribution of the error over time, we present graphs shown in Figure \ref{errorsC4} and Figure \ref{errorsC4_2}. The graphs are constructed for $C=10^4$ and contain the results for 3-7 oversampling layers.  
We get a significant error jump at the beginning with a coarse grid $20 \times $20. The error jump is small when we use a coarse grid of $40 \times 40$. The error behavior became smooth after the beginning. We observe no error jumps throughout the process, keeping the error at the same level. We notice errors increasing throughout the simulation for three oversampling layers, which explains why we should employ four or more oversampling layers.

We demonstrated that the presented NLMC algorithm gives highly accurate calculation results. Even in highly heterogeneous domains, the approach demonstrated great precision. The method's accuracy is practically independent of the nonlinear part of the equation and the type of heterogeneity. We can improve accuracy by adding more oversampling layers. We observe that the method depends on the coarse grid's size. The approach is more accurate on a finer coarse grid.
Furthermore, the method significantly reduces the original system's size. In our experiments, the coarse grid system size is significantly smaller than the fine grid system size ($562 \times 562$ or $1930 \times 1930$ for the coarse grid and $41733 \times 41733$ for the fine grid). From this point, we can see that multiscale methods save computational resources, which is the primary advantage of multiscale methods over traditional mathematical modeling approaches. 


\section{Conclusion}

This paper presented a Non-Local Multi-Continuum method (NLMC) algorithm for the time-dependent Darcy-Forchheimer model in a fractured heterogeneous domain.   The fine grid approximation was constructed using the Finite volume method with a lower dimensional embedded fracture model. Model and methods formulations were given for two-dimensional cases. We completed a numerical experiment including two test cases with different heterogeneous properties. The numerical results are obtained for varying numbers of oversampling layers. We investigated nonlinearity's impact by varying the coefficient $\beta$. The numerical experiment showed that the proposed approach had provided accurate results without significant influence of the nonlinear part of the flow. The performed numerical experiment showed good accuracy of the method. We conclude that the Non-Local Multi-Continuum technique performed well in modeling the Darcy-Forchheimer model in the fractured medium.

\section{Acknowledgements}

D. Spiridonov work is supported the grant of Russian Science Foundation No. 21-71-00061(https://rscf.ru/en/project/21-71-00061/) and the Russian government project Science and Universities (project No. FSRG-2021-0015) aimed at supporting junior laboratories.

\bibliographystyle{unsrt}
\bibliography{lit}

\end{document}